# The Age-Structured Chemostat with Substrate Dynamics as a Control System

**Iasson Karafyllis[*], Dionysis Theodosis[*] and Miroslav Krstic[**]**

[*]Dept. of Mathematics, National Technical University of Athens, Zografou Campus, 15780, Athens, Greece,
emails: iasonkar@central.ntua.gr , dtheodp@central.ntua.gr

[**]Dept. of Mechanical and Aerospace Eng., University of California, San Diego, La Jolla, CA 92093-0411, U.S.A., email: krstic@ucsd.edu

## Abstract

In this work we study an age-structured chemostat model with a renewal boundary condition and a coupled substrate equation. The model is nonlinear and consists of a hyperbolic partial differential equation and an ordinary differential equation with nonlinear, nonlocal terms appearing both in the ordinary differential equation and the boundary condition. Both differential equations contain a non-negative control input, while the states of the model are required to be positive. Under an appropriate weak solution framework, we determine the state space and the input space for this model. We prove global existence and uniqueness of solutions for all admissible initial conditions and all allowable control inputs. To this purpose we employ a combination of Banach's fixed-point theorem with implicit solution formulas and useful solution estimates. Finally, we show that the age-structured chemostat model gives a well-defined control system on a metric space.

## 1. Introduction

A chemostat is a continuous-flow bioreactor which has been widely used for maintaining continuous microbial growth. In a chemostat, fresh medium enters and culture leaves at the same rate, keeping the working volume constant. The chemostat has found numerous applications in practical systems like wastewater treatment, biomass, biofuel, and pharmaceuticals production. The chemostat has a central role in mathematical biology, where it has been studied both as a dynamical system and as a control system. On the dynamical systems side, numerous studies examine the stability properties of equilibria and long-term behavior of the chemostat such as existence of periodic solution under periodic substrate for both single-species and multiple competitive organisms (see for instance [1], [4], [6], [9], [24], [25], [27], [33]). From a control viewpoint, the chemostat has been studied through a variety of feedback designs and methods that leverage the dilution rate as control input. These methods include LQ-based strategies, adaptive feedback schemes under uncertainties, stabilization under measurement and actuation delays, control Lyapunov function, as well as dynamic feedback law designs (see for instance [2], [7], [8], [12], [14], [20], [21]).

A key topic in mathematical modeling for biology, medicine, demography, and economics is the analysis of structured population models that encode distributions over age, size, and sex within a single population. Age-structured models account for the distribution of individuals across their life cycle where, instead of treating the population as homogeneous, they represent how processes like growth, division, and death depend on age, and they include a boundary condition that generates new individuals when older ones reproduce or divide (see for instance [3], [10], [29], [32]). A classic example of age-structured model is the well-known as McKendrick-von Foerster ([11], [22]), which is a first order hyperbolic Partial Differential Equation (PDE) with a nonlocal boundary condition. In age-structured models, ergodic theorems are sometimes used to characterize long-time behavior, see [15], [16], [17], [19].

The mathematical model of the simple chemostat has several limitations. As quoted by Pilyugin and Waltman in [23], "Following the accumulation of experimental data, it became evident that the simple model requires modification. Specifically, the simple model failed to explain the observed oscillatory behavior in the chemostat." Age-structured models are natural extensions of the standard chemostat: they keep the same continuous-flow idea but also include age as a second variable capturing how growth depends on the age of the population (see [13], [19], [26], [32]). This additional structure enriches the model's qualitative behavior by capturing oscillations that simple Ordinary Differential Equation (ODE) models fail to reproduce, see [30], [31]. Thus, age-structured formulations therefore provide a more realistic and dynamically richer description of continuous bioreactors. However, when the age structure of the microbial population is coupled to substrate dynamics, the resulting PDE–ODE system introduces nonlinear and nonlocal terms that destroy the necessary structure for the application of the ergodic theorem (as was done in [15], [16], [17], [19]).

In this paper we deal with the age-structured chemostat model with substrate dynamics proposed in [30], [31] which consists of a first-order PDE coupled with an ODE describing the substrate rate of change and a nonlinear, non-local boundary condition. This is the first paper that studies the age-structured chemostat with substrate dynamics as a control system under very mild and physically reasonable assumptions. It should be noted that the stabilization problem of the age-structured chemostat was also studied in [20]. However, [20] required some demanding structural assumptions that led to a simpler finite-dimensional system. The contributions of the paper are as follows:

1) The first contribution of the paper is the clarification of the notion of solution for the age-structured chemostat model with substrate dynamics. The proposed notion of a weak-solution retains sufficient regularity for the study of the qualitative properties of the chemostat.

2) The second contribution of the paper is the well-posedness result. We show that the

corresponding initial-boundary value problem is well-posed under mild assumptions which have important physical interpretations (e.g., non-negativity and boundedness of the birth modulus, positivity of the state, etc.). Moreover, we show that the problem admits a global solution for all physically meaningful initial conditions and inputs.

3) The third -and most important from a control-theoretic point of view- is the clarification of the state space, the input space and the formulation of the age-structured chemostat as a well-defined control system.

Having performed the above contributions, future researchers can exploit the proposed formulation and obtain stability and stabilization results or controllability and optimal control results. In other words, our work provides the initial step for a detailed study of the age-structured chemostat with substrate dynamics.

The paper is outlined as follows. Section 2 presents the age-structured chemostat model, its state space and an appropriate definition of a weak solution. Section 3 contains the main results that include global existence and uniqueness of solutions of the age-structured chemostat model, as well as the chemostat as a control system. Finally, Section 4 is devoted to the proofs of the main results, while some concluding remarks are given in Section 5.

**Notation** Throughout this paper, we adopt the following notation.

* $\mathbb{R}_+ := [0, +\infty)$. For a vector $x \in \mathbb{R}^n$, $|x|$ denotes its Euclidean norm.

* Let $A \subseteq \mathbb{R}^n$ be an open set and let $B \subseteq \mathbb{R}^n$ be a set that satisfies $A \subseteq B \subseteq cl(A)$, where $cl(A)$ is the closure of $A$. By $C^0(B; \Omega)$, we denote the class of continuous functions on $B$, which take values in $\Omega \subseteq \mathbb{R}^m$. By $C^k(B; \Omega)$, where $k \geq 1$ is an integer, we denote the class of functions on $B \subseteq \mathbb{R}^n$, which take values in $\Omega \subseteq \mathbb{R}^m$ and have continuous derivatives of order $k$. In other words, the functions of class $C^k(B; \Omega)$ are the functions which have continuous derivatives of order $k$ in $A = int(B)$ that can be continued continuously to all points in $\partial A \cap B$. When $\Omega = \mathbb{R}$ then we write $C^0(B)$ or $C^k(B)$.

* Let $A \subseteq \mathbb{R}^n$ be an open set and let $\Omega \subseteq \mathbb{R}^m$ be a non-empty set. By $L^p(A; \Omega)$ with $p \geq 1$ we denote the equivalence class of measurable functions $f : A \to \Omega$ for which $\|f\|_p = \left( \int_A |f(x)|^p \, dx \right)^{1/p} < +\infty$. By $L^\infty(A; \Omega)$ we denote the equivalence class of measurable functions $f : A \to \Omega$ for which $\|f\|_\infty = \sup_{x \in A} (|f(x)|) < +\infty$ where $\sup_{x \in A} (|f(x)|)$ is the essential supremum. When $\Omega = \mathbb{R}^m$ we simply write $L^p(A)$. When $B \subseteq \mathbb{R}^n$ is not open but has non-empty interior, $L^p(B; \Omega)$ and $L^\infty(B; \Omega)$ mean $L^p(A; \Omega)$ and $L^\infty(A; \Omega)$, respectively, with $A = int(B)$. By $L^\infty_{loc}(\mathbb{R}_+; \Omega)$ we denote

the equivalence class of measurable functions $f: \mathbb{R}_+ \to \Omega$ with $f \in L^\infty([0,T]; \Omega)$ for every $T > 0$.

* Let $I \subseteq \mathbb{R}_+$ and let $f: I \times \mathbb{R}_+ \to \mathbb{R}_+$ be a given function. We use the notation $f[t]$ to denote the profile at certain $t \in I$, i.e., $(f[t])(a) = f(t,a)$ for all $a \geq 0$.

## 2. The Age-Structured Chemostat Model

Consider the age-structured chemostat model

$$\frac{\partial f}{\partial t}(t,a) + \frac{\partial f}{\partial a}(t,a) = -\left(\beta(a) + D(t)\right) f(t,a), \qquad (2.1)$$

$$f(t,0) = \mu(S(t)) \int_0^{+\infty} k(a) f(t,a) da, \qquad (2.2)$$

$$\dot{S}(t) = D(t)\left(S_{in} - S(t)\right) - \mu(S(t)) \int_0^{+\infty} q(a) f(t,a) da, \qquad (2.3)$$

where $f(t,a) > 0$ is the distribution function of the microbial population in the chemostat at time $t \geq 0$ and age $a \geq 0$, $S(t) > 0$ is the limiting substrate concentration, $S_{in} > 0$ is the inlet concentration of the substrate, $D(t) > 0$ is the dilution rate, $\mu(S)$ is the specific growth rate function, $\beta(a)$ is the mortality rate and $k(a), q(a)$ are functions that determine the birth of new cells and the substrate consumption of the microbial population, respectively. All functions $\mu, \beta, k, q : \mathbb{R}_+ \to \mathbb{R}_+$ are assumed to be bounded, $C^0(\mathbb{R}_+)$ functions with $\mu \in C^1(\mathbb{R}_+)$, $\mu(0) = 0$, $\mu(S) > 0$ for $S > 0$ and $\int_0^{+\infty} k(a) da > 0$, $\int_0^{+\infty} q(a) da > 0$.

Clearly, system (2.1), (2.2), (2.3) is a complicated nonlinear model that consists of a hyperbolic PDE (known as McKendrick–von Foerster equation, see [11], [22]) with a non-local boundary condition and an ODE. Boundedness of $\mu, \beta, k, q : \mathbb{R}_+ \to \mathbb{R}_+$ reflects finite physiological capacities, such as no infinite growth, mortality, birth and substrate consumption, while $\int_0^{+\infty} k(a) da > 0$ and $\int_0^{+\infty} q(a) da > 0$ imply that there is at least some age range with non-zero reproduction and some age range where substrate is actually consumed. Common choices on the growth rate are the Monod kinetics $\mu(S) = \mu_{\max} S / (K_S + S)$, $\mu_{\max}, K_S > 0$ and the Haldane kinetics $\mu(S) = \mu_{\max} S / (K_P + S + S^2 K_I^{-1})$, $K_P, K_I > 0$ (see [27], [9]).

Define:

$$X = \left\{ (f,S) \in L^1(\mathbb{R}_+) \times (0, S_{in}) : \begin{array}{c} f \in C^0(\mathbb{R}_+;(0,+\infty)), \lim_{a \to +\infty}(f(a)) = 0 \\ f(0) = \mu(S) \int_0^{+\infty} k(a) f(a) da \end{array} \right\} \quad (2.4)$$

We consider $X$ to be a metric space with metric given by the formula for all $(f,S), (\bar{f}, \bar{S}) \in X$:

$$d\left((f,S), (\bar{f}, \bar{S})\right) = \|f - \bar{f}\|_1 + |S - \bar{S}| \quad (2.5)$$

The results that are given below show that the metric space $X$ with metric given by (2.5) is the state space of model (2.1), (2.2), (2.3).

We next provide the notion of solution that is appropriate for system (2.1), (2.2), (2.3).

**Definition 1:** *Let $D \in L^\infty_{loc}(\mathbb{R}_+; \mathbb{R}_+)$, $(f_0, S_0) \in X$ and $T > 0$ be given. We say that a continuous mapping $(f,S):[0,T] \to X$ is a weak solution on $[0,T]$ with input $D$ of the initial-boundary value problem (2.1), (2.2), (2.3) with initial condition*

$$f[0] = f_0, \, S(0) = S_0 \quad (2.6)$$

*if the following properties are valid:*
*i) (2.6) holds,*
*ii) $f \in C^0([0,T] \times \mathbb{R}_+;(0,+\infty))$ and $S:[0,T] \to (0, S_{in})$ is absolutely continuous,*
*iii) (2.2) holds for all $t \in [0,T]$ and (2.3) holds for $t \in [0,T]$ a.e.*
*iv) the following equation holds for all $\varphi \in C^1([0,T] \times \mathbb{R}_+) \cap L^\infty([0,T] \times \mathbb{R}_+)$ with $\left( \frac{\partial \varphi}{\partial a} + \frac{\partial \varphi}{\partial t} \right) \in L^\infty([0,T] \times \mathbb{R}_+)$ and $t \in [0,T]$:*

$$\int_0^{+\infty} f_0(a) \varphi(0,a) da + \int_0^t f(s,0) \varphi(s,0) ds = \int_0^{+\infty} f(t,a) \varphi(t,a) da$$
$$+ \int_0^t \int_0^{+\infty} \left( (\beta(a) + D(s)) \varphi(s,a) - \frac{\partial \varphi}{\partial a}(s,a) - \frac{\partial \varphi}{\partial s}(s,a) \right) f(s,a) da ds \quad (2.7)$$

*A weak solution on $[0,T]$ with input $D$ of the initial-boundary value problem (2.1), (2.2), (2.3), (2.6) $(f,S):[0,T] \to X$ is called a classical solution on $[0,T]$ of the initial-boundary value problem (2.1), (2.2), (2.3), (2.6) if $f \in C^1([0,T] \times \mathbb{R}_+)$, $S \in C^1([0,T];(0, S_{in}))$ and (2.1), (2.3) hold for all $(t,a) \in [0,T] \times \mathbb{R}_+$.*

*We say that the initial-boundary value problem (2.1), (2.2), (2.3), (2.6) admits a global weak (classical) solution with input $D$ if for every $T > 0$ there exists a weak (classical) solution on $[0,T]$ with input $D$ of the initial-boundary value problem (2.1), (2.2), (2.3), (2.6).*

**Remark 1:** It becomes clear from Definition 1 that for every $T > 0$, $\varphi \in C^1(\mathbb{R}_+) \cap L^\infty(\mathbb{R}_+)$ with $\varphi' \in L^\infty(\mathbb{R}_+)$ and for every weak solution on $[0,T]$ with input $D$ of the initial-boundary value problem (2.1), (2.2), (2.3), (2.6) the following equation holds for $t \in [0,T]$ a.e.:

$$f(t,0)\varphi(0) + \int_0^{+\infty} \left(\varphi'(a) - (\beta(a) + D(t))\varphi(a)\right) f(t,a)\,da = \frac{d}{dt}\int_0^{+\infty} f(t,a)\varphi(a)\,da \qquad (2.8)$$

## 3. Main Results

Having clarified the notion of the solution, we are in a position to present the main results of this paper. The following theorem establishes global existence of solutions of the initial value problem (2.1), (2.2), (2.3), (2.6).

**Theorem 1:** *For every $(f_0, S_0) \in X$ and $D \in L^\infty_{loc}(\mathbb{R}_+; \mathbb{R}_+)$, there exists a global weak solution with input $D$ of the initial-boundary value problem (2.1), (2.2), (2.3), (2.6).*

The proof of Theorem 1 is provided in Section 4. Its proof employs a combination of Banach's fixed-point theorem with implicit solution formulas and some useful solution estimates.

The following result shows continuous dependence on the initial conditions.

**Theorem 2:** *Suppose that $(f, S)$ is a weak solution on $[0,T]$ of the initial-boundary value problem (2.1), (2.2), (2.3), (2.6) with input $D \in L^\infty_{loc}(\mathbb{R}_+; \mathbb{R}_+)$ for certain $T > 0$. Then there exists a constant $\chi > 0$ that depends only on $\sup_{s \in [0,T]}(D(s))$ and $\max_{s \in [0,T]}(\|f[s]\|_1)$ such that the following estimate holds for every weak solution $(\tilde{f}, \tilde{S})$ on $[0,T]$ of the initial-boundary value problem (2.1), (2.2), (2.3), with $(\tilde{f}[0], \tilde{S}(0)) = (\tilde{f}_0, \tilde{S}_0) \in X$ and input $\tilde{D} \in L^\infty_{loc}(\mathbb{R}_+; \mathbb{R}_+)$ with $\tilde{D}(t) = D(t)$ for $t \in [0,T]$ a.e.:*

$$\|f[t] - \tilde{f}[t]\|_1 + |S(t) - \tilde{S}(t)| \leq \exp(\chi t)\left(\|f_0 - \tilde{f}_0\|_1 + |S_0 - \tilde{S}_0|\right), \text{ for all } t \in [0,T] \quad (3.1)$$

Having established continuous dependence on the initial conditions, uniqueness of the weak solution for a given input $D \in L^\infty_{loc}(\mathbb{R}_+; \mathbb{R}_+)$ and a given initial condition $(f_0, S_0) \in X$ is a consequence of Theorem 2.

**Corollary 1:** *For any $(f_0, S_0) \in X$, $D \in L^\infty_{loc}(\mathbb{R}_+; \mathbb{R}_+)$, there exists a unique weak solution $(f, S)$ on $[0, T]$, $T > 0$ with input $D$ of the initial-boundary value problem (2.1), (2.2), (2.3), (2.6).*

Therefore, the age-structured chemostat initial-boundary value problem (2.1), (2.2), (2.3), (2.6) is a *well-posed problem in the sense of Hadamard* (see page 155 in [18]).

Theorem 1 allows us to associate to each input $D \in L^\infty_{loc}(\mathbb{R}_+; \mathbb{R}_+)$ and $(f_0, S_0) \in X$ a well-defined solution $(f[t], S[t]) \in X$. We define for all $D \in L^\infty_{loc}(\mathbb{R}_+; \mathbb{R}_+)$, $(f_0, S_0) \in X$ and $t \geq 0$:

$$\phi(t, (f_0, S_0); D) = (f[t], S[t]) \tag{3.2}$$

Theorem 1, Theorem 2 and Definition 1 allow us to guarantee that the mapping $\phi$ defined by (3.2) is a mapping

$$\phi : \mathbb{R}_+ \times X \times L^\infty_{loc}(\mathbb{R}_+; \mathbb{R}_+) \to X$$

that satisfies the identity property

$$\phi(0, (f_0, S_0); D) = (f_0, S_0), \text{ for all } D \in L^\infty_{loc}(\mathbb{R}_+; \mathbb{R}_+), (f_0, S_0) \in X \tag{3.3}$$

and the causality property

$$\phi(t, (f_0, S_0); D) = \phi(t, (f_0, S_0); \tilde{D})$$
$$\text{for all } t > 0, (f_0, S_0) \in X \text{ and } D, \tilde{D} \in L^\infty_{loc}(\mathbb{R}_+; \mathbb{R}_+)$$
$$\text{with } \tilde{D}(s) = D(s) \text{ for } s \in [0, t] \text{ a.e.} \tag{3.4}$$

Moreover, $\phi(t, (f_0, S_0); D)$ is continuous with respect to $t \geq 0$ and $(f_0, S_0) \in X$.

The classical semigroup property

$$\phi(t + \tau, (f_0, S_0); D) = \phi(t, \phi(\tau, (f_0, S_0); D); \delta_\tau D),$$
$$\text{for all } t, \tau \geq 0, D \in L^\infty_{loc}(\mathbb{R}_+; \mathbb{R}_+), (f_0, S_0) \in X \tag{3.5}$$

where $\delta_\tau D$ is the shifted input $(\delta_\tau D)(s) = D(\tau + s)$, is a consequence of the following technical proposition.

**Proposition 1:** *Let $T, \tau > 0$ and $D \in L^\infty_{loc}(\mathbb{R}_+; \mathbb{R}_+)$. Assume that $(f, S)$ is a weak solution on $[0, T]$ with input $D$ of the initial-boundary value problem (2.1), (2.2), (2.3), (2.6) and assume that $(\bar{f}, \bar{S})$ is a weak solution on $[0, \tau]$ with input $\bar{D}(s) := D(T + s)$ of the initial-boundary value problem (2.1), (2.2), (2.3), with $\bar{f}[0] = f[T]$, $\bar{S}(0) = S(T)$. Define*

$$\left(\hat{f}[t], \hat{S}(t)\right) = \begin{cases} \left(f[t], S(t)\right) & 0 \leq t \leq T \\ \left(\bar{f}[t-T], \bar{S}(t-T)\right) & T < t \leq T+\tau \end{cases} \quad (3.6)$$

*Then $(\hat{f}, \hat{S})$ is a weak solution on $[0, T+\tau]$ with input $D$ of the initial-boundary value problem (2.1), (2.2), (2.3), (2.6).*

Thus, following the terminology in [28], the age-structured chemostat model (2.1), (2.2), (2.3), defines a *forward complete time-invariant control system with state space $X$ and input space $L^\infty_{loc}(\mathbb{R}_+; \mathbb{R}_+)$.*

## 4. Proofs

For reader's convenience, we first present an outline of this section. First, we prove a local existence theorem and we establish continuous dependence of solutions on initial conditions. Then we exploit some estimates of the solutions, and we extend the solutions globally.

We start with the two following technical lemmas. Their proofs are provided in the Appendix.

**Lemma 1:** *Let $T > 0$, $x \in C^0([0,T])$, $D \in L^\infty([0,T])$, $f_0 \in C^0(\mathbb{R}_+) \cap L^1(\mathbb{R}_+)$ and $\beta \in C^0(\mathbb{R}_+) \cap L^\infty(\mathbb{R}_+)$ with $f_0(0) = x(0)$. The function $f \in C^0([0,T] \times \mathbb{R}_+)$ with $\sup_{t \in [0,T]}\left(\|f[t]\|_1\right) < +\infty$ defined by*

$$f(t,a) = \begin{cases} f_0(a-t)\exp\left(-\int_0^t D(s)ds - \int_{a-t}^a \beta(s)ds\right) & \text{for} \quad 0 \leq t \leq a \\ x(t-a)\exp\left(-\int_{t-a}^t D(s)ds - \int_0^a \beta(s)ds\right) & \text{for} \quad t > a \geq 0 \end{cases} \quad (4.1)$$

*is the unique function in $\left\{u \in C^0([0,T] \times \mathbb{R}_+): \sup_{t \in [0,T]}\left(\|u[t]\|_1\right) < +\infty\right\}$ that satisfies the following equation for all $\varphi \in C^1([0,T] \times \mathbb{R}_+) \cap L^\infty([0,T] \times \mathbb{R}_+)$ with $\left(\frac{\partial \varphi}{\partial a} + \frac{\partial \varphi}{\partial t}\right) \in L^\infty([0,T] \times \mathbb{R}_+)$ and $t \in [0,T]$:*

$$\int_0^{+\infty} f_0(a)\varphi(0,a)da + \int_0^t x(s)\varphi(s,0)ds = \int_0^{+\infty} f(t,a)\varphi(t,a)da$$
$$+ \int_0^t \int_0^{+\infty} \left( (\beta(a)+D(s))\varphi(s,a) - \frac{\partial \varphi}{\partial a}(s,a) - \frac{\partial \varphi}{\partial s}(s,a) \right) f(s,a)dads \quad (4.2)$$

**Lemma 2:** *Let $D \in L^\infty_{loc}(\mathbb{R}_+;\mathbb{R}_+)$, $(f_0,S_0) \in X$ and $T > 0$ be given. Let $(f,S):[0,T] \to X$ be a weak solution on $[0,T]$ with input $D$ of the initial-boundary value problem (2.1), (2.2), (2.3), (2.6). Let $\Gamma > 0$ be a constant for which the inequality $\mu(S) \leq \Gamma S$ holds for all $S \in [0,S_{in}]$ and define $M = \max_{S \in [0,S_{in}]}(\mu(S))$. Then the following estimates hold for all $t \in [0,T]$:*

$$\|f[t]\|_1 \leq \exp(M \|k\|_\infty t) \|f_0\|_1 \quad (4.3)$$

$$S_0 \exp\left( -\Gamma \|q\|_\infty \|f_0\|_1 \frac{\exp(M\|k\|_\infty t) - 1}{M \|k\|_\infty} \right) \leq S(t)$$

$$S(t) \leq S_{in}\left( 1 - \exp\left(-\int_0^t D(s)ds\right) \right) + S_0 \exp\left(-\int_0^t D(s)ds\right) \quad (4.4)$$

We next provide a local existence result.

**Theorem 3:** *There exists a continuous function $\Delta:(0,S_{in}) \times \mathbb{R}_+ \to (0,1]$ with the following property: for every $(f_0,S_0) \in X$, $D \in L^\infty_{loc}(\mathbb{R}_+;\mathbb{R}_+)$ and $\delta \in \left( 0, \Delta\left(S_0, \|f_0\|_1 + \sup_{s \in [0,1]}(D(s))\right) \right]$ there exists a weak solution on $[0,\delta]$ with input $D$ of the initial-boundary value problem (2.1), (2.2), (2.3), (2.6).*

**Proof:** Let $D \in L^\infty_{loc}(\mathbb{R}_+;\mathbb{R}_+)$ and $(f_0,S_0) \in X$ be given. Define:

$$\bar{g}(t) := \int_t^{+\infty} k(a) f_0(a-t) \exp\left(-\int_{a-t}^a \beta(s)ds\right) da$$

$$\bar{h}(t) := \int_t^{+\infty} q(a) f_0(a-t) \exp\left(-\int_{a-t}^a \beta(s)ds\right) da \quad (4.5)$$

$$b(t) := \exp\left(-\int_0^t D(s)ds\right)$$

The functions $\bar{g}, \bar{h}$ are continuous. To see this, notice that definition (4.5) gives

$$\bar{g}(t) = \int_0^{+\infty} p(t,u)\,du$$

$$p(t,u) := k(u+t)f_0(u)\exp\left(-\int_u^{u+t}\beta(s)ds\right), u \geq 0$$

Continuity of $k, \beta$ and boundedness of $\beta$, gives that for any fixed $t \geq 0$ and for any sequence $\{t_n : n = 0,1,2,...\}$ with $\lim_{n\to+\infty}(t_n) = t$, and $u \geq 0$

$$\lim_{n\to+\infty}(p(t_n,u)) = p(t,u)$$

Since $\beta(s) \geq 0$ for all $s \geq 0$ and $c = \sup_{s\geq 0}(k(s)) < +\infty$, we have for all $n \geq 0$ and for all $u \geq 0$

$$|p(t_n,u)| \leq c|f_0(u)|$$

Taking into account that $f_0 \in L^1((0,+\infty))$ and the Dominated Convergence Theorem we get that

$$\lim_{n\to+\infty}(\bar{g}(t_n)) = \lim_{n\to+\infty}\left(\int_0^{+\infty} p(t_n,u)du\right) = \left(\int_0^{+\infty}\lim_{n\to+\infty}(p(t_n,u))du\right) = \left(\int_0^{+\infty} p(t,u)du\right) = \bar{g}(t)$$

Since the latter holds for any sequence $\{t_n : n = 0,1,2,...\}$ and any $t \geq 0$, we conclude that $\bar{g} \in C^0(\mathbb{R}_+)$. Completely analogous is the proof of continuity of $\bar{h}$.

Let time $\delta \in (0,1]$ which is to be selected in an appropriate way below.

In what follows, we consider the functions $\bar{g}, \bar{h}$, given by (4.5) to be defined on $[0,1]$. Thus, when we write $\|\bar{g}\|_\infty$ and $\|\bar{h}\|_\infty$ we mean $\max_{t\in[0,1]}(|\bar{g}(t)|)$ and $\max_{t\in[0,1]}(|\bar{h}(t)|)$, respectively.

We also define the functions:

$$\tilde{k}(a) := k(a)\exp\left(-\int_0^a \beta(s)ds\right), \quad \tilde{q}(a) := q(a)\exp\left(-\int_0^a \beta(s)ds\right) \quad (4.6)$$

Consider the Banach space $H = C^0([0,\delta]; \mathbb{R}^2)$ with norm $\|(y,z)\|_H = \|y\|_\infty + \|z\|_\infty = \max_{t\in[0,\delta]}(|y(t)|) + \max_{t\in[0,\delta]}(|z(t)|)$. Define $R := \frac{1}{2}\min(S_0, S_{in} - S_0) > 0$ and notice that due to the fact that $S_0 \in (0, S_{in})$ (recall definition (2.4) and the fact that $(f_0, S_0) \in X$) we get that for every $(y,z) \in H$ with $\|(y,z)\|_H \leq R$ we obtain $0 < \frac{S_0}{2} \leq S_0 + z(t) \leq \frac{S_{in} + S_0}{2} < S_{in}$ for all $t \in [0,\delta]$. Define

$$B_R = \left\{(y,z) \in H : \left\|(y,z)\right\|_H \leq R\right\} \tag{4.7}$$

which is a closed subset of $H$ and the operator

$$T = (T_1, T_2) : B_R \to H \tag{4.8}$$

by the following formulas

$$(T_1(y,z))(t) = \int_0^t \tilde{k}(a)\mu(z(t-a)+S_0)(y(t-a)+\bar{g}(t-a))da \tag{4.9}$$

$$(T_2(y,z))(t) = (S_{in} - S_0)(1-b(t)) - b(t)\int_0^t \mu(z(\tau)+S_0)\bar{h}(\tau)d\tau$$
$$-b(t)\int_0^t \mu(z(\tau)+S_0)\left(\int_0^\tau \tilde{q}(a)\mu(z(\tau-a)+S_0)(y(\tau-a)+\bar{g}(\tau-a))da\right)d\tau \tag{4.10}$$

We show next, by appropriately choosing $\delta > 0$, that $T(B_R) \subset B_R$.

Since $\mu(s) \geq 0$ for all $s \geq 0$ we get from (4.9) for all $(y,z) \in B_R$:

$$\left\|T_1(y,z)\right\|_\infty \leq M \left\|\tilde{k}\right\|_\infty \left(R + \left\|\bar{g}\right\|_\infty\right)\delta \tag{4.11}$$

where

$$M = \sup_{S \geq 0}(\mu(S)) \tag{4.12}$$

Since $D(t) \geq 0$ for $t \in [0,\delta]$ a.e. and $\delta \leq 1$, we get from (4.5) for all $t \in [0,\delta]$

$$b(t) \leq 1 \text{ and } 1-b(t) \leq \int_0^t D(s)ds \leq \delta \sup_{s \in [0,1]}(D(s)) \tag{4.13}$$

Using (4.12), (4.13) and the fact that $\mu(s) \geq 0$ for all $s \geq 0$, we obtain from (4.10) for all $(y,z) \in B_R$:

$$\left\|T_2(y,z)\right\|_\infty \leq \left((S_{in} - S_0)\sup_{s \in [0,1]}(D(s)) + M\left\|\bar{h}\right\|_\infty + M^2\left\|\tilde{q}\right\|_\infty\left(R + \left\|\bar{g}\right\|_\infty\right)\frac{\delta}{2}\right)\delta \tag{4.14}$$

It follows from (4.11) and (4.14) and the fact that $R := \frac{1}{2}\min(S_0, S_{in} - S_0) \leq S_{in}$ that $\left\|T_1(y,z)\right\|_\infty + \left\|T_2(y,z)\right\|_\infty \leq R$ for all $\delta \in (0,1]$ with

$$\delta \leq \min\left(1, \frac{\min(S_0, S_{in} - S_0)/2}{S_{in} \sup_{s \in [0,1]}(D(s)) + M \|\bar{h}\|_\infty + M(\|\tilde{k}\|_\infty + M\|\tilde{q}\|_\infty)(S_{in} + \|\bar{g}\|_\infty)}\right) \quad (4.15)$$

Let arbitrary $(y_1, z_1), (y_2, z_2) \in B_R$ be given. Using (4.12), (4.14), the fact that $\mu(s) \geq 0$ for all $s \geq 0$ and by adding and subtracting terms in (4.9), (4.10), we obtain for all $t \in [0, \delta]$:

$$\begin{aligned}
&|T_1(y_1, z_1)(t) - T_1(y_2, z_2)(t)| \\
&\leq \int_0^t |\tilde{k}(a)||\mu(S_0 + z_1(t-a)) - \mu(S_0 + z_2(t-a))||y_2(t-a) + \bar{g}(t-a)|da \\
&\quad + \int_0^t |\tilde{k}(a)|\mu(S_0 + z_1(t-a))|y_1(t-a) - y_2(t-a)|da \\
&\leq \|\tilde{k}\|_\infty L_\mu (R + \|\bar{g}\|_\infty) \delta \|z_1 - z_2\|_\infty + M\|\tilde{k}\|_\infty \delta \|y_1 - y_2\|_\infty
\end{aligned} \quad (4.16)$$

$$\begin{aligned}
&\exp\left(\int_0^t D(s)ds\right)|T_2(y_1, z_1)(t) - T_2(y_2, z_2)(t)| \\
&\leq \int_0^t |\mu(z_1(\tau) + S_0) - \mu(z_2(\tau) + S_0)||\bar{h}(\tau)|d\tau \\
&\quad + \int_0^t |\mu(z_1(\tau) + S_0) - \mu(z_2(\tau) + S_0)| \int_0^\tau \tilde{q}(a)\mu(z_1(\tau-a) + S_0)|y_1(\tau-a) + \bar{g}(\tau-a)|dad\tau \\
&\quad + \int_0^t \mu(z_2(\tau) + S_0) \int_0^\tau \tilde{q}(a)|\mu(z_1(\tau-a) + S_0) - \mu(z_2(\tau-a) + S_0)||y_1(\tau-a) + \bar{g}(\tau-a)|dad\tau \\
&\quad + \int_0^t \mu(z_2(\tau) + S_0) \int_0^\tau \tilde{q}(a)\mu(z_2(\tau-a) + S_0)|y_1(\tau-a) - y_2(\tau-a)|dad\tau \\
&\leq L_\mu \|\bar{h}\|_\infty \delta \|z_1 - z_2\|_\infty + M\|\tilde{q}\|_\infty L_\mu \|z_1 - z_2\|_\infty (R + \|\bar{g}\|_\infty)\delta^2 + M^2 \|\tilde{q}\|_\infty \|y_1 - y_2\|_\infty \frac{\delta^2}{2}
\end{aligned}$$
(4.17)

where $L_\mu = \max\{|\mu'(s)| : 0 \leq s \leq S_{in}\}$. Therefore, we get from (4.16), (4.17) and the fact that $\int_0^t D(s)ds \geq 0$ for all $t \geq 0$:

$$\begin{aligned}
&\|T_1(y_1, z_1) - T_1(y_2, z_2)\|_\infty + \|T_2(y_1, z_1) - T_2(y_2, z_2)\|_\infty \\
&\leq M\left(\|\tilde{k}\|_\infty + M\|\tilde{q}\|_\infty \delta\right) \delta \|y_1 - y_2\|_\infty \\
&\quad + L_\mu \left(\|\bar{h}\|_\infty + \left(M\|\tilde{q}\|_\infty \delta + \|\tilde{k}\|_\infty\right)(R + \|\bar{g}\|_\infty)\right) \delta \|z_1 - z_2\|_\infty
\end{aligned} \quad (4.18)$$

The facts that $\|\bar{h}\|_\infty \leq \|q\|_\infty \|f_0\|_1, \|\bar{g}\|_\infty \leq \|k\|_\infty \|f_0\|_1, \|\tilde{k}\|_\infty \leq \|k\|_\infty, \|\tilde{q}\|_\infty \leq \|q\|_\infty$ (direct consequences of definitions (4.6), (4.5) and the fact that $\beta$ is non-negative), in

conjunction with (4.15), (4.18) show that the operator $T = (T_1, T_2) : B_R \to H$ defined by (4.9), (4.10) is a contraction for every $\delta \in \left[ 0, \Delta\left( S_0, \|f_0\|_1 + \sup_{s \in [0,1]} (D(s)) \right) \right]$ with $\Delta(s, r) := \dfrac{\min(s, S_{in} - s)}{2KS_{in}(r+1)}$ for all $s > 0, r \geq 0$, where $K > 0$ is a sufficiently large constant independent of $D \in L^\infty_{loc}(\mathbb{R}_+; \mathbb{R}_+)$ and $(f_0, S_0) \in X$, for example

$$K = 1 + \left( \frac{M}{S_{in}} + L_\mu \right)\left( \|q\|_\infty + \left( \|k\|_\infty + M \|q\|_\infty \right)\|k\|_\infty \right) + \left( M + L_\mu S_{in} \right)\left( M \|q\|_\infty + \|k\|_\infty \right)$$

Thus, by Banach's fixed point theorem there exists $(y, z) \in B_R$ such that the following equations hold for all $t \in [0, \delta]$:

$$y(t) = \int_0^t \tilde{k}(a) \mu(z(t-a) + S_0)(y(t-a) + \bar{g}(t-a)) da \tag{4.19}$$

$$\begin{aligned}
z(t) &= (S_{in} - S_0)(1 - b(t)) - b(t) \int_0^t \mu(z(\tau) + S_0) \bar{h}(\tau) d\tau \\
&\quad - b(t) \int_0^t \mu(z(\tau) + S_0) \left( \int_0^\tau \tilde{q}(a) \mu(z(\tau - a) + S_0)(y(\tau - a) + \bar{g}(\tau - a)) da \right) d\tau
\end{aligned} \tag{4.20}$$

Define for $t \in [0, \delta]$ and $a \geq 0$:

$$\begin{aligned}
S(t) &= z(t) + S_0 \\
x(t) &= \mu(S(t)) b(t) (\bar{g}(t) + y(t))
\end{aligned} \tag{4.21}$$

$$f(t, a) = \begin{cases} f_0(a - t) \exp\left( -\int_0^t D(s) ds - \int_{a-t}^a \beta(s) ds \right) & \text{for } 0 \leq t \leq a \\ x(t - a) \exp\left( -\int_{t-a}^t D(s) ds - \int_0^a \beta(s) ds \right) & \text{for } t > a \geq 0 \end{cases} \tag{4.22}$$

We show next that $(f, S)$ is a weak solution on $[0, \delta]$ with input $D$ of the initial-boundary value problem (2.1), (2.2), (2.3), (2.6).

We proceed by showing some facts.

Fact 0: Equation (2.6) holds.

Fact 0 is a consequence of equations (4.5), (4.20), (4.21), (4.22).

<u>Fact 1:</u> $S:[0,\delta]\to(0,S_{in})$ is absolutely continuous.

Proof of Fact 1: Equations (4.5), (4.20), (4.21) show that $S:[0,\delta]\to\mathbb{R}$ is absolutely continuous. Furthermore, the facts that $S_0\in(0,S_{in})$ (recall definition (2.4) and the fact that $(f_0,S_0)\in X$), $(y,z)\in B_R$ and $R=\frac{1}{2}\min(S_0,S_{in}-S_0)$ in conjunction with (4.7) and (4.21) imply that $0<\frac{S_0}{2}\le S(t)\le\frac{S_{in}+S_0}{2}<S_{in}$ for all $t\in[0,\delta]$. Therefore, $S:[0,\delta]\to(0,S_{in})$.

<u>Fact 2:</u> $x\in C^0([0,\delta])$ and $f\in C^0([0,\delta]\times\mathbb{R}_+)$ with

$$f(t,0)=x(t) \text{ and } \lim_{a\to+\infty}(f(t,a))=0, \text{ for all } t\in[0,\delta] \qquad (4.23)$$

Proof of Fact 2: Since $\bar{g}$ is continuous, it follows from (4.21) that $x\in C^0([0,\delta])$. Moreover, equations (4.5), (4.19), (4.20), (4.21) imply that $x(0)=\mu(S_0)\int_0^{+\infty}k(a)f_0(a)da$. Definition (2.4) and the fact that $(f_0,S_0)\in X$ imply that $x(0)=f_0(0)$. Therefore, it follows from (4.22) and the fact that $(f_0,S_0)\in X$ (which implies that $f_0\in C^0(\mathbb{R}_+;(0,+\infty))$; recall (2.4)) that $f\in C^0([0,\delta]\times\mathbb{R}_+)$. Moreover, (4.22) and the fact that $\lim_{a\to+\infty}(f_0(a))=0$ (a consequence of the fact that $(f_0,S_0)\in X$; recall (2.4)) guarantees that (4.23) holds.

<u>Fact 3:</u> $f[t]\in L^1(\mathbb{R}_+)\cap C^0(\mathbb{R}_+;(0,+\infty))$, for all $t\in[0,\delta]$ with $\sup_{t\in[0,\delta]}(\|f[t]\|_1)<+\infty$.

Proof of Fact 3: Using (4.12), (4.13), (4.7), (4.5), (4.21) and the fact that $\beta(a)\ge0$ for all $a\ge0$ we get:
$$\|x\|_\infty\le M(\|k\|_\infty\|f_0\|_1+R) \qquad (4.24)$$

Using (4.22), (4.24) and the facts that $\beta(a)\ge0$ for all $a\ge0$, $D(t)\ge0$ for $t\ge0$ a.e., we get for $t\in[0,\delta]$:

$$\|f[t]\|_1 = \int_0^{+\infty} |f(t,a)|\,da = \int_0^t |f(t,a)|\,da + \int_t^{+\infty} |f(t,a)|\,da$$

$$= \int_0^t |x(t-a)|\exp\left(-\int_{t-a}^t D(s)\,ds - \int_0^a \beta(s)\,ds\right) da$$

$$+ \int_t^{+\infty} |f_0(a-t)|\exp\left(-\int_0^t D(s)\,ds - \int_{a-t}^a \beta(s)\,ds\right) da$$

$$\leq \int_0^t |x(t-a)|\,da + \int_t^{+\infty} |f_0(a-t)|\,da \leq \|x\|_\infty t + \|f_0\|_1$$

$$\leq M\left(\|k\|_\infty \|f_0\|_1 + R\right)\delta + \|f_0\|_1$$

Therefore, for every $t \in [0,\delta]$ it holds that $f[t] \in L^1(\mathbb{R}_+) \cap C^0(\mathbb{R}_+;(0,+\infty))$ with $\sup_{t\in[0,\delta]}\left(\|f[t]\|_1\right) \leq M\left(\|k\|_\infty \|f_0\|_1 + R\right)\delta + \|f_0\|_1 < +\infty$.

<u>Fact 4</u>: Equation (2.2) holds for all $t \in [0,\delta]$.

Proof of Fact 4: Using (4.5) and (4.19), (4.21), (4.23) we get for all $t \in [0,\delta]$:

$$\begin{aligned}f(t,0) = x(t) &= \mu(S(t))\exp\left(-\int_0^t D(s)\,ds\right)\int_t^{+\infty} k(a)f_0(a-t)\exp\left(-\int_{a-t}^a \beta(s)\,ds\right) da \\ &+ \mu(S(t))\exp\left(-\int_0^t D(s)\,ds\right)\int_0^t \tilde{k}(a)\mu(S(t-a))\left(y(t-a)+\overline{g}(t-a)\right) da\end{aligned} \quad (4.25)$$

Using (4.6), (4.5), (4.21), (4.22) we get for all $t \in [0,\delta]$:

$$\int_0^{+\infty} k(a)f(t,a)\,da = \int_0^t k(a)x(t-a)\exp\left(-\int_{t-a}^t D(s)\,ds - \int_0^a \beta(s)\,ds\right) da$$

$$+ \int_t^{+\infty} k(a)f_0(a-t)\exp\left(-\int_0^t D(s)\,ds - \int_{a-t}^a \beta(s)\,ds\right) da$$

$$= \int_0^t \tilde{k}(a)x(t-a)\exp\left(-\int_{t-a}^t D(s)\,ds\right) da$$

$$+ \exp\left(-\int_0^t D(s)\,ds\right)\int_t^{+\infty} k(a)f_0(a-t)\exp\left(-\int_{a-t}^a \beta(s)\,ds\right) da$$

$$= \exp\left(-\int_0^t D(s)\,ds\right)\int_0^t \tilde{k}(a)\mu(S(t-a))\left(\overline{g}(t-a)+y(t-a)\right) da$$

$$+ \exp\left(-\int_0^t D(s)\,ds\right)\int_t^{+\infty} k(a)f_0(a-t)\exp\left(-\int_{a-t}^a \beta(s)\,ds\right) da$$

The above equation and (4.25) show that equation (2.2) holds for all $t \in [0, \delta]$.

<u>Fact 5</u>: $f(t,a) > 0$ for all $t \in [0, \delta]$, $a \geq 0$.

Proof of Fact 5: Since $f_0(a) > 0$ for all $a \geq 0$ (a consequence of the fact that $(f_0, S_0) \in X$; recall (2.4)), equation (4.22) shows that it suffices to show that $x(t) > 0$ for all $t \in [0, \delta]$.

Since $\mu(S) > 0$ for $S > 0$ and since $S \in C^0([0, \delta]; (0, S_{in}))$ (recall Fact 1), equations (4.5) and (4.21) show that it suffices to show that $\bar{g}(t) + y(t) > 0$ for all $t \in [0, \delta]$.

We notice that equations (4.5) and (4.19) show that $\bar{g}(0) + y(0) = \int_0^{+\infty} k(a) f_0(a) da > 0$ which is a consequence of the facts that $f_0(a) > 0$, $k(a) \geq 0$ for all $a \geq 0$ with $\int_0^{+\infty} k(a) da > 0$.

The proof of the fact that $\bar{g}(t) + y(t) > 0$ for all $t \in [0, \delta]$ is made by contradiction. Suppose that there exists $t \in [0, \delta]$ with $\bar{g}(t) + y(t) \leq 0$. Then there exists $T \in (0, t]$ with $\bar{g}(T) + y(T) = 0$ and $\bar{g}(s) + y(s) > 0$ for all $s \in [0, T)$. Indeed, the set $\{\tau \in [0, T] : \bar{g}(\tau) + y(\tau) = 0\}$ is bounded and non-empty (a consequence of the facts that $\bar{g}(0) + y(0) > 0$, $\bar{g}(t) + y(t) \leq 0$ and Bolzano's theorem). Setting $T = \inf \{\tau \in [0, T] : \bar{g}(\tau) + y(\tau) = 0\}$, by continuity we get that $\bar{g}(T) + y(T) = 0$ and $\bar{g}(s) + y(s) > 0$ for all $s \in [0, T)$.

Equations (4.5) and (4.19) show that

$$0 = \bar{g}(T) + y(T) = \int_T^{+\infty} k(a) f_0(a-T) \exp\left(-\int_{a-T}^a \beta(s) ds\right) da \\ + \int_0^T k(a) \exp\left(-\int_0^a \beta(s) ds\right) \mu(S(T-a))(y(T-a) + \bar{g}(T-a)) da \quad (4.26)$$

Since $f_0(a) > 0$, $k(a) \geq 0$ for all $a \geq 0$, $\mu(S) > 0$ for $S > 0$, $S \in C^0([0, \delta]; (0, S_{in}))$ and $\bar{g}(s) + y(s) > 0$ for all $s \in [0, T)$, equation (4.26) shows that

$$\int_T^{+\infty} k(a) f_0(a-T) \exp\left(-\int_{a-T}^a \beta(s) ds\right) da = 0$$

$$\int_0^T k(a) \exp\left(-\int_0^a \beta(s) ds\right) \mu(S(T-a))(y(T-a) + \bar{g}(T-a)) da = 0$$

The above equations can only hold if $k(a) = 0$ for all $a \geq 0$; a contradiction with the fact that $\int_0^{+\infty} k(a)da > 0$.

A consequence of all the previous facts and definition (2.4) is the following fact.

<u>Fact 6:</u> $(f[t], S(t)) \in X$ for all $t \in [0, \delta]$.

The following fact is a direct consequence of definition (4.22), Fact 2, Fact 3 and Lemma 1.

<u>Fact 7:</u> Equation (2.7) holds for all $\varphi \in C^1([0,\delta] \times \mathbb{R}_+) \cap L^\infty([0,\delta] \times \mathbb{R}_+)$ with $\left(\frac{\partial \varphi}{\partial a} + \frac{\partial \varphi}{\partial t}\right) \in L^\infty([0,\delta] \times \mathbb{R}_+)$ and $t \in [0, \delta]$.

We next show the following continuity result.

<u>Fact 8:</u> The mapping $[0, \delta] \ni t \to f[t] \in L^1(\mathbb{R}_+)$ is continuous.

Proof of Fact 8: Let arbitrary $\varepsilon > 0$ and $t_0 \in [0, \delta]$ be given.

We note that there exists $A \geq 1$ such that $\int_A^{+\infty} |f(s,a)| da < \varepsilon / 4$ for all $s \in [0, \delta]$. Indeed, since $f_0 \in L^1(\mathbb{R}_+)$ there exists $\bar{A} \geq 0$ such that $\int_{\bar{A}}^{+\infty} |f_0(a)| da < \varepsilon / 4$. Taking $A = 1 + \bar{A}$ and using the facts that $\delta \leq 1$ (which implies that $s \leq a$ for all $s \in [0, \delta]$ and $a \geq A$), $D, \beta$ are non-negative functions, we get from (4.22) for $s \in [0, \delta]$:

$$\int_A^{+\infty} |f(s,a)| da = \int_A^{+\infty} |f_0(a-s)| \exp\left(-\int_0^s D(l)dl - \int_{a-s}^a \beta(l)dl\right) da$$

$$\leq \int_A^{+\infty} |f_0(a-s)| da = \int_{1+\bar{A}-s}^{+\infty} |f_0(r)| dr \leq \int_{\bar{A}}^{+\infty} |f_0(r)| dr < \varepsilon / 4$$

Thus, we get for all $t \in [0, \delta]$:

$$\|f[t] - f[t_0]\|_1 = \int_0^A |f(t,a) - f(t_0,a)| da + \int_A^{+\infty} |f(t,a) - f(t_0,a)| da$$

$$\leq \int_0^A |f(t,a) - f(t_0,a)| da + \int_A^{+\infty} |f(t,a)| da + \int_A^{+\infty} |f(t_0,a)| da \quad (4.27)$$

$$< \int_0^A |f(t,a) - f(t_0,a)| da + \varepsilon / 2$$

Fact 2 (and consequently uniform continuity on the compact set $[0,\delta]\times[0,A]$) implies the existence of $h>0$ such that $|f(t,a)-f(t_0,a)|<\dfrac{\varepsilon}{2A}$ for all $t\in[0,\delta]$ with $|t-t_0|<h$. Inequality (4.27) implies that $\|f[t]-f[t_0]\|_1<\varepsilon$ for all $t\in[0,\delta]$ with $|t-t_0|<h$.

Clearly, a direct consequence of Fact 8, Fact 1 and definition (2.5) is the fact that the mapping $(f,S):[0,T]\to X$ is continuous.

The only thing that remains to be proved is the validity of equation (2.3) for $t\in[0,\delta]$ a.e..

Using (4.20), (4.21), and the definitions of $\bar{h},b$ in (4.5) we obtain

$$S(t) = S_{in} - b(t)(S_{in}-S_0)$$
$$-b(t)\int_0^t \mu(S(\tau))\left(\int_\tau^{+\infty} q(a)f_0(a-\tau)\exp\left(-\int_{a-\tau}^a \beta(s)ds\right)da\right)d\tau \qquad (4.28)$$
$$-b(t)\int_0^t \mu(S(\tau))\left(\int_0^\tau \tilde{q}(a)x(\tau-a)b^{-1}(\tau-a)da\right)d\tau$$

From (4.28) and the definitions of $\tilde{q},b$ in (4.6) and (4.5), respectively, it follows that

$$S(t) = S_{in} - b(t)(S_{in}-S_0)$$
$$-b(t)\int_0^t \mu(S(\tau))\left(\int_\tau^{+\infty} q(a)f_0(a-\tau)\exp\left(-\int_{a-\tau}^a \beta(s)ds\right)b^{-1}(\tau)b(\tau)da\right)d\tau$$
$$-b(t)\int_0^t \mu(S(\tau))\left(\int_0^\tau q(a)x(\tau-a)\exp\left(-\int_0^a \beta(s)ds\right)b^{-1}(\tau-a)b^{-1}(\tau)b(\tau)da\right)d\tau$$
$$= S_{in} - b(t)(S_{in}-S_0)$$
$$-b(t)\int_0^t \mu(S(\tau))\left(\int_\tau^{+\infty} q(a)f_0(a-\tau)\exp\left(-\int_{a-\tau}^a \beta(s)ds\right)b^{-1}(\tau)b(\tau)da\right)d\tau$$
$$-b(t)\int_0^t \mu(S(\tau))\left(\int_0^\tau q(a)x(\tau-a)\exp\left(-\int_0^a \beta(s)ds\right)\exp\left(-\int_{\tau-a}^\tau D(s)ds\right)b^{-1}(\tau)da\right)d\tau$$

Taking into account (4.22) and the definition of $b$ in (4.5), the previous equality gives that

$$S(t) = S_{in} - b(t)(S_{in}-S_0) - b(t)\int_0^t \mu(S(\tau))\left(\int_0^\infty q(a)f(\tau,a)da\right)b^{-1}(\tau)d\tau \qquad (4.29)$$

Due to Fact 8, and boundedness of $q$, it follows that $t \mapsto \int_0^\infty q(a) f(t,a) da$ is continuous on $[0,\delta]$. Then, the validity of (2.3) for $t \in [0,\delta]$ a.e. follows directly from (4.5) and (4.29).

The proof is complete. ◁

We continue with the proof of Theorem 2.

**Proof of Theorem 2:** Let arbitrary weak solution $(\tilde{f}, \tilde{S})$ on $[0,T]$ of the initial-boundary value problem (2.1), (2.2), (2.3), with initial condition $(\tilde{f}[0], \tilde{S}(0)) = (\tilde{f}_0, \tilde{S}_0) \in X$ and input $\tilde{D} \in L^\infty_{loc}(\mathbb{R}_+; \mathbb{R}_+)$ with $\tilde{D}(t) = D(t)$ for $t \in [0,T]$ a.e. be given. By virtue of Lemma 1, we have for $t \in [0,T]$ that

$$f(t,a) = \begin{cases} f_0(a-t) \exp\left(-\int_0^t D(s)ds - \int_{a-t}^a \beta(s)ds\right) & \text{for } 0 \le t \le a \\ x(t-a) \exp\left(-\int_{t-a}^t D(s)ds - \int_0^a \beta(s)ds\right) & \text{for } t > a \ge 0 \end{cases} \quad (4.30)$$

$$\tilde{f}(t,a) = \begin{cases} \tilde{f}_0(a-t) \exp\left(-\int_0^t D(s)ds - \int_{a-t}^a \beta(s)ds\right) & \text{for } 0 \le t \le a \\ \tilde{x}(t-a) \exp\left(-\int_{t-a}^t D(s)ds - \int_0^a \beta(s)ds\right) & \text{for } t > a \ge 0 \end{cases} \quad (4.31)$$

where $x(t) = f(t,0)$ and $\tilde{x}(t) = \tilde{f}(t,0)$. For notational convenience, we define

$$\begin{aligned} \Delta f[t] &:= f[t] - \tilde{f}[t] \\ \Delta S(t) &:= S(t) - \tilde{S}(t) \\ \Delta x(t) &:= x(t) - \tilde{x}(t) \end{aligned} \quad (4.32)$$

Since $D(t) \ge 0$ for $t \ge 0$ a.e. and since $\beta(s) \ge 0$ for $s \ge 0$, we get from (4.30), (4.31) and (4.32) for $t \in [0,T]$

$$\begin{aligned} \|f[t] - \tilde{f}[t]\|_1 &= \|\Delta f[t]\|_1 = \int_0^t |\Delta f(t,a)| da + \int_t^{+\infty} |\Delta f(t,a)| da \le \int_0^t |\Delta x(s)| ds \\ &+ \int_0^t |f_0(a-t) - \tilde{f}_0(a-t)| da \le \|f_0 - \tilde{f}_0\|_1 + \int_0^t |\Delta x(s)| ds \end{aligned} \quad (4.33)$$

Using the boundary condition (2.2) for $x(t) = f(t,0)$ and $\tilde{x}(t) = \tilde{f}(t,0)$, and (4.32) we have for $t \in [0,T]$ that

$$|\Delta x(t)| \leq |\mu(S(t)) - \mu(\tilde{S}(t))| \int_0^{+\infty} k(a) f(t,a) da + \mu(\tilde{S}(t)) \left| \int_0^{+\infty} k(a) \Delta f(t,a) da \right| \quad (4.34)$$

Since $\mu \in C^1(\mathbb{R}_+)$, we have that

$$|\mu(v) - \mu(r)| \leq L_\mu |v - r| \quad (4.35)$$

for all $v, r \in [0, S_{in}]$, where $L_\mu := \max\{|\mu'(s)| : 0 \leq s \leq S_{in}\}$. From (4.35), (4.34), and definition $M = \max_{s \in [0, S_{in}]} (\mu(s))$ we get that

$$|\Delta x(t)| \leq L_\mu \|k\|_\infty \max_{s \in [0,T]} (\|f[s]\|_1) |\Delta S(t)| + M \|k\|_\infty \|\Delta f[t]\|_1 \quad (4.36)$$

where $\max_{s \in [0,T]} (\|f[s]\|_1) < +\infty$ (a consequence of continuity of $(f, S) : [0,T] \to X$ and definition (2.5) which implies continuity of $f : [0,T] \to L^1(\mathbb{R}_+)$).

Since both $S(t)$ and $\tilde{S}(t)$ satisfy (2.3) and since $\tilde{D}(t) = D(t)$ for $t \in [0,T]$ a.e., we get for $t \in [0,T]$ a.e.

$$\frac{d}{dt}(S(t) - \tilde{S}(t)) = \frac{d}{dt}(\Delta S(t)) = -D(t) \Delta S(t)$$
$$-\mu(\tilde{S}(t)) \int_0^{+\infty} q(a) (f(t,a) - \tilde{f}(t,a)) da - (\mu(S(t)) - \mu(\tilde{S}(t))) \int_0^{+\infty} q(a) f(t,a) da$$

Integrating the above equation, we get for all $t \in [0,T]$ that

$$\Delta S(t) = S_0 - \tilde{S}_0 - \int_0^t D(s) \Delta S(s) ds - \int_0^t \mu(\tilde{S}(s)) \int_0^{+\infty} q(a) (f(s,a) - \tilde{f}(s,a)) da\, ds$$
$$- \int_0^t (\mu(S(s)) - \mu(\tilde{S}(s))) \int_0^{+\infty} q(a) f(s,a) da\, ds \quad (4.37)$$

Using the (4.37), together with (4.32), (4.35), and $\max_{t \in [0,T]} (\|f[t]\|_1) < +\infty$ we have that

$$|\Delta S(t)| \leq |S_0 - \tilde{S}_0| + \left( \sup_{s \in [0,T]} (D(s)) + L_\mu \|q\|_\infty \max_{s \in [0,T]} (\|f[s]\|_1) \right) \int_0^t |\Delta S(s)| ds$$
$$+ M \|q\|_\infty \int_0^t \|\Delta f[s]\|_1 ds \quad (4.38)$$

Combining (4.33), (4.36), and (4.38) we get for $t \in [0,T]$ that

$$\|\Delta f[t]\|_1 + |\Delta S(t)| \le \|f_0 - \tilde{f}_0\|_1 + |S_0 - \tilde{S}_0| + \int_0^t |\Delta x(s)| ds$$

$$+ \left( \sup_{s\in[0,T]} (D(s)) + L_\mu \|q\|_\infty \max_{s\in[0,T]} (\|f[s]\|_1) \right) \int_0^t |\Delta S(s)| ds + M \|q\|_\infty \int_0^t \|\Delta f[s]\|_1 ds$$

$$\le \|f_0 - \tilde{f}_0\|_1 + |S_0 - \tilde{S}_0| + \chi \int_0^t \left( \|\Delta f[s]\|_1 + |\Delta S(s)| \right) ds$$

where $\chi := \sup_{s\in[0,T]} (D(s)) + \left( \max_{s\in[0,T]} (\|f[s]\|_1) L_\mu + M \right) (\|k\|_\infty + \|q\|_\infty)$. A direct application of Gronwall's inequality gives (3.1). The proof is complete. ◁

Proposition 1 shows how we can extend the solution of a given initial-boundary value problem. Its proof is provided below.

**Proof of Proposition 1:** We show that for all $\varphi \in C^1([0,T+\tau]\times\mathbb{R}_+) \cap L^\infty([0,T+\tau]\times\mathbb{R}_+)$ with $\left( \frac{\partial \varphi}{\partial a} + \frac{\partial \varphi}{\partial t} \right) \in L^\infty([0,T+\tau]\times\mathbb{R}_+)$ and $t \in [0, T+\tau]$ the following equation holds:

$$\int_0^{+\infty} f_0(a)\varphi(0,a)da + \int_0^t \hat{f}(s,0)\varphi(s,0)ds = \int_0^{+\infty} \hat{f}(t,a)\varphi(t,a)da \quad (4.39)$$

$$+ \int_0^t \int_0^{+\infty} \left( (\beta(a) + D(s))\varphi(s,a) - \frac{\partial \varphi}{\partial a}(s,a) - \frac{\partial \varphi}{\partial s}(s,a) \right) \hat{f}(s,a) da ds$$

Let arbitrary $\varphi \in C^1([0,T+\tau]\times\mathbb{R}_+) \cap L^\infty([0,T+\tau]\times\mathbb{R}_+)$ with $\left( \frac{\partial \varphi}{\partial a} + \frac{\partial \varphi}{\partial t} \right) \in L^\infty([0,T+\tau]\times\mathbb{R}_+)$ be given. Then $\varphi \in C^1([0,T]\times\mathbb{R}_+) \cap L^\infty([0,T]\times\mathbb{R}_+)$ with $\left( \frac{\partial \varphi}{\partial a} + \frac{\partial \varphi}{\partial t} \right) \in L^\infty([0,T+\tau]\times\mathbb{R}_+)$ and since $(f, S)$ is a weak solution on $[0,T]$ with input $D$ of the initial-boundary value problem (2.1), (2.2), (2.3), (2.6) it follows (from Definition 1) that (2.7) holds for $t \in [0,T]$. Therefore, definition (3.6) guarantees that (4.39) holds for $t \in [0,T]$.

Since $\varphi \in C^1([0,T+\tau]\times\mathbb{R}_+) \cap L^\infty([0,T+\tau]\times\mathbb{R}_+)$ with $\left( \frac{\partial \varphi}{\partial a} + \frac{\partial \varphi}{\partial t} \right) \in L^\infty([0,T+\tau]\times\mathbb{R}_+)$, it follows that the function

$$\bar{\varphi}[s] = \varphi[T+s], \text{ for } s \in [0,\tau] \quad (4.40)$$

is a function of class $\bar{\varphi} \in C^1([0,\tau] \times \mathbb{R}_+) \cap L^\infty([0,\tau] \times \mathbb{R}_+)$ with $\left(\frac{\partial \bar{\varphi}}{\partial a} + \frac{\partial \bar{\varphi}}{\partial t}\right) \in L^\infty([0,\tau] \times \mathbb{R}_+)$. Since $(\bar{f}, \bar{S})$ is a weak solution on $[0,\tau]$ with input $\bar{D}(s) := D(T+s)$ of the initial-boundary value problem (2.1), (2.2), (2.3), with $\bar{f}[0] = f[T]$, $\bar{S}(0) = S(T)$, it follows (from Definition 1) that the following equation holds for all $t \in [0,\tau]$:

$$\int_0^{+\infty} f(T,a)\bar{\varphi}(0,a)da + \int_0^t \bar{f}(s,0)\bar{\varphi}(s,0)ds = \int_0^{+\infty} \bar{f}(t,a)\bar{\varphi}(t,a)da \\ + \int_0^t \int_0^{+\infty} \left((\beta(a) + D(T+s))\bar{\varphi}(s,a) - \frac{\partial \bar{\varphi}}{\partial a}(s,a) - \frac{\partial \bar{\varphi}}{\partial s}(s,a)\right)\bar{f}(s,a)dads \quad (4.41)$$

Using definitions (3.6), (4.40), we obtain from (4.41) for all $t \in (0,\tau]$:

$$\int_0^{+\infty} f(T,a)\varphi(T,a)da + \int_T^{T+t} \hat{f}(s,0)\varphi(s,0)ds = \int_0^{+\infty} \hat{f}(T+t,a)\varphi(T+t,a)da \\ + \int_T^{T+t} \int_0^{+\infty} \left((\beta(a) + D(s))\varphi(s,a) - \frac{\partial \varphi}{\partial a}(s,a) - \frac{\partial \varphi}{\partial s}(s,a)\right)\hat{f}(s,a)dads \quad (4.42)$$

Exploiting (2.7) for $t = T$ and definition (3.6) we get:

$$\int_0^{+\infty} f_0(a)\varphi(0,a)da + \int_0^T \hat{f}(s,0)\varphi(s,0)ds = \int_0^{+\infty} f(T,a)\varphi(T,a)da \\ + \int_0^T \int_0^{+\infty} \left((\beta(a) + D(s))\varphi(s,a) - \frac{\partial \varphi}{\partial a}(s,a) - \frac{\partial \varphi}{\partial s}(s,a)\right)\hat{f}(s,a)dads \quad (4.43)$$

Combining (4.42) and (4.43) we obtain (4.39) for all $t \in (T, T+\tau]$.

All the rest requirements of Definition 1 are direct consequences of definition (3.6). The proof is complete. ◁

Finally, we can provide the proof of Theorem 1.

**Proof of Theorem 1:** Let arbitrary $(f_0, S_0) \in X$, $D \in L^\infty_{loc}(\mathbb{R}_+; \mathbb{R}_+)$ be given. Define the set

$$J := \left\{ \delta > 0 : \begin{array}{l} \text{a weak solution with input } D \text{ of} \\ (1), (2), (3), (7) \text{ exists on } [0, \delta] \end{array} \right\} \quad (4.44)$$

and let

$$T_{\max} = \sup(J) \quad (4.45)$$

Notice that by virtue of Theorem 3, $J \neq \emptyset$ and thus, $T_{\max} \in (0, +\infty]$. We next show by contradiction that $T_{\max} = +\infty$.

Assume that $T_{\max} < +\infty$. Define,

$$a := S_0 \exp\left(-\Gamma \|q\|_\infty \|f_0\|_1 \frac{\exp(M\|k\|_\infty T_{\max}) - 1}{M\|k\|_\infty}\right) > 0$$

$$b := S_{in} - (S_{in} - S_0)\exp\left(-\int_0^{T_{\max}} D(s)ds\right) < S_{in}$$

and

$$F^* := \exp(M\|k\|_\infty T_{\max})\|f_0\|_1$$

Since $D \in L^\infty_{loc}(\mathbb{R}_+; \mathbb{R}_+)$, it follows that $\sup_{t \in [0, T_{\max}+1]} (D(t)) < +\infty$ and thus the set

$$K := [a,b] \times \left[0, F^* + \sup_{t \in [0, T_{\max}+1]}(D(t))\right] \subset (0, S_{in}) \times \mathbb{R}_+$$

is compact. Consider now the continuous function $\Delta : (0, S_{in}) \times \mathbb{R}_+ \to (0,1]$ provided by Theorem 3. Continuity of $\Delta$ on the compact set $K$ gives

$$\delta^* := \min_{(s, \ell) \in K}(\Delta(s, \ell)) > 0 \tag{4.46}$$

Let $\{t_n \in J : n = 0,1,2,...\}$ be a non-decreasing sequence with $\lim_{n \to +\infty}(t_n) = T_{\max}$. Let $N > 0$ be such that $t_N > T_{\max} - \delta^*/2$ (this is possible since $\lim_{n \to +\infty}(t_n) = T_{\max}$). Since $t_N \in J$, it follows from definition (4.44) that a weak solution $(f[t], S(t))$ of the initial-boundary value problem (2.1), (2.2), (2.3), (2.6) exists on $[0, t_N]$. Then by virtue of Lemma 2, for all $t \in [0, t_n]$ it holds that

$$\|f[t]\|_1 \leq F^* \\ 0 < a \leq S(t) \leq b < S_{in} \tag{4.47}$$

Moreover, from Theorem 3, a weak solution $(\tilde{f}, \tilde{S})$ of the initial-boundary value problem (2.1), (2.2), (2.3), with $(\tilde{f}[0], \tilde{S}(0)) = (f[t_N], S(t_N)) \in X$ and input $D_N(s) := D(t_N + s)$, $s \geq 0$ exists on $[0, \bar{\delta}]$ for every $\bar{\delta} \in \left(0, \Delta\left(S(t_N), \|f[t_N]\|_1 + \sup_{s \in [t_N, t_N+1]}(D(s))\right)\right)$. Due to (4.47) and definition of $K$, it follows that $\left(S(t_N), \|f[t_N]\|_1 + \sup_{s \in [t_N, t_N+1]} (D(s))\right) \in K$. Therefore, from (4.46) we obtain that

$$\delta^* \leq \Delta\left(S(t_N), \|f[t_N]\|_1 + \max_{s \in [t_N, t_N+1]} (D(s))\right)$$

We pick $\bar{\delta} = \delta^*$. From Proposition 1, it follows that

$$(\hat{f}[t], \hat{S}(t)) := \begin{cases} (f[t], S(t)), & t \in [0, t_N] \\ \left(\tilde{f}[t - t_N], S(t - t_N)\right), & t \in [t_N, t_N + \delta^*] \end{cases}$$

is a weak solution on $[0, t_N + \delta^*]$ with input $D$ of the initial-boundary value problem (2.1), (2.2), (2.3), (2.6), where $t_N + \delta^* > T_{max}$ (recall that $t_N > T_{max} - \delta^*/2$). The latter contradicts definition (4.45) of $T_{max}$. Thus, $T_{max} = +\infty$. The proof is complete. ◁

## 5. Conclusions

We have studied the well-posedness of an age-structured chemostat model with a nonlocal (renewal) boundary condition and a coupled substrate equation. Under an appropriate weak solution framework, we have determined the state space and the input space for this model, and we have proved global existence and uniqueness of solutions for all admissible initial conditions and all allowable control inputs. Our formulation of the age-structured chemostat as a well-posed control system, opens the door to a systematic study of stability and stabilization, and enables subsequent advances in controllability and optimal control.

# Appendix

**Proof of Lemma 1:** The fact that the function $f$ defined by (4.1) is of class $C^0([0,T]\times\mathbb{R}_+)$ is a consequence of the facts that $x \in C^0([0,T])$, $D \in L^\infty([0,T])$ and $f_0, \beta \in C^0(\mathbb{R}_+)$ with $f_0(0) = x(0)$. We also have from (4.1) and the facts that $x \in C^0([0,T])$, $D \in L^\infty([0,T])$, $f_0 \in C^0(\mathbb{R}_+) \cap L^1(\mathbb{R}_+)$ and $\beta \in C^0(\mathbb{R}_+) \cap L^\infty(\mathbb{R}_+)$ for every $t \in [0,T]$:

$$\|f[t]\|_1 = \int_0^t |f(t,a)|da + \int_t^{+\infty} |f(t,a)|da$$

$$= \int_0^t |x(t-a)| \exp\left(-\int_{t-a}^t D(l)dl - \int_0^a \beta(l)dl\right) da$$

$$+ \int_t^{+\infty} |f_0(a-t)| \exp\left(-\int_0^t D(l)dl - \int_{a-t}^a \beta(l)dl\right) da$$

$$\leq \int_0^t |x(t-a)| \exp\left((\|D\|_\infty + \|\beta\|_\infty)a\right) da$$

$$+ \int_t^{+\infty} |f_0(a-t)| \exp\left((\|D\|_\infty + \|\beta\|_\infty)t\right) da$$

$$\leq \exp\left((\|D\|_\infty + \|\beta\|_\infty)T\right)\left(\|x\|_\infty T + \|f_0\|_1\right)$$

Thus, the function $f$ defined by (4.1) satisfies

$$\sup_{t \in [0,T]} \left(\|f[t]\|_1\right) \leq \exp\left((\|D\|_\infty + \|\beta\|_\infty)T\right)\left(\|x\|_\infty T + \|f_0\|_1\right) < +\infty$$

We next show that the function $f$ defined by (4.1) satisfies equation (4.2) for all $\varphi \in C^1([0,T]\times\mathbb{R}_+) \cap L^\infty([0,T]\times\mathbb{R}_+)$ with $\left(\frac{\partial\varphi}{\partial a} + \frac{\partial\varphi}{\partial t}\right) \in L^\infty([0,T]\times\mathbb{R}_+)$ and $t \in [0,T]$.

Let (arbitrary) $t \in [0,T]$ and $\varphi \in C^1([0,T]\times\mathbb{R}_+) \cap L^\infty([0,T]\times\mathbb{R}_+)$ with $\left(\frac{\partial\varphi}{\partial a} + \frac{\partial\varphi}{\partial t}\right) \in L^\infty([0,T]\times\mathbb{R}_+)$ be given. We get:

$$\int_0^t \int_0^{+\infty} \left((\beta(a) + D(s))\varphi(s,a) - \frac{\partial\varphi}{\partial a}(s,a) - \frac{\partial\varphi}{\partial s}(s,a)\right) f(s,a) da\, ds$$
$$+ \int_0^{+\infty} f(t,a)\varphi(t,a) da = I_1 + I_2 \quad \text{(A1)}$$

where

$$I_1 := \int_0^t f(t,a)\varphi(t,a)da + \int_0^t\int_0^s \Omega(s,a)f(s,a)dads$$

$$I_2 := \int_t^{+\infty} f(t,a)\varphi(t,a)da + \int_0^t\int_s^{+\infty} \Omega(s,a)f(s,a)dads \qquad (A2)$$

$$\Omega(s,a) := (\beta(a)+D(s))\varphi(s,a) - \frac{\partial\varphi}{\partial a}(s,a) - \frac{\partial\varphi}{\partial s}(s,a) \qquad (A3)$$

Definition (4.1) in conjunction with definitions (A2) imply that:

$$
\begin{aligned}
I_1 &= \int_0^t x(t-a)\exp\left(-\int_{t-a}^t D(s)ds - \int_0^a \beta(s)ds\right)\varphi(t,a)da \\
&+ \int_0^t\int_0^s \Omega(s,a)x(s-a)\exp\left(-\int_{s-a}^s D(l)dl - \int_0^a \beta(l)dl\right)dads \\
I_2 &= \int_t^{+\infty} f_0(a-t)\exp\left(-\int_0^t D(s)ds - \int_{a-t}^a \beta(s)ds\right)\varphi(t,a)da \\
&+ \int_0^t\int_s^{+\infty} \Omega(s,a)f_0(a-s)\exp\left(-\int_0^s D(l)dl - \int_{a-s}^a \beta(l)dl\right)dads
\end{aligned} \qquad (A4)
$$

Using (A4), Fubini's theorem and definition (A3) we get:

$$
\begin{aligned}
I_2 &= \int_0^t\int_0^{+\infty} \Omega(s,s+r)f_0(r)\exp\left(-\int_0^s D(l)dl - \int_r^{s+r}\beta(l)dl\right)drds \\
&+ \int_0^{+\infty} f_0(r)\exp\left(-\int_0^t D(s)ds - \int_r^{t+r}\beta(s)ds\right)\varphi(t,t+r)dr \\
&= \int_0^{+\infty} f_0(r)\int_0^t \Omega(s,s+r)\exp\left(-\int_0^s D(l)dl - \int_r^{s+r}\beta(l)dl\right)dsdr \\
&+ \int_0^{+\infty} f_0(r)\exp\left(-\int_0^t D(s)ds - \int_r^{t+r}\beta(s)ds\right)\varphi(t,t+r)dr \\
&= -\int_0^{+\infty} f_0(r)\int_0^t \frac{\partial}{\partial s}\left(\exp\left(-\int_0^s D(l)dl - \int_r^{s+r}\beta(l)dl\right)\varphi(s,s+r)\right)dsdr \\
&+ \int_0^{+\infty} f_0(r)\exp\left(-\int_0^t D(s)ds - \int_r^{t+r}\beta(s)ds\right)\varphi(t,t+r)dr = \int_0^{+\infty} f_0(r)\varphi(0,r)dr
\end{aligned}
$$

$$I_1 = \int_0^t \int_0^s \Omega(s, s-r) x(r) \exp\left(-\int_r^s D(l)dl - \int_0^{s-r} \beta(l)dl\right) dr\,ds$$

$$+ \int_0^t x(r) \exp\left(-\int_r^t D(s)ds - \int_0^{t-r} \beta(s)ds\right) \varphi(t, t-r) dr$$

$$= \int_0^t \int_r^t \Omega(s, s-r) x(r) \exp\left(-\int_r^s D(l)dl - \int_0^{s-r} \beta(l)dl\right) ds\,dr$$

$$+ \int_0^t x(r) \exp\left(-\int_r^t D(s)ds - \int_0^{t-r} \beta(s)ds\right) \varphi(t, t-r) dr$$

$$= -\int_0^t x(r) \int_r^t \frac{\partial}{\partial s}\left(\varphi(s, s-r) \exp\left(-\int_r^s D(l)dl - \int_0^{s-r} \beta(l)dl\right)\right) ds\,dr$$

$$+ \int_0^t x(r) \exp\left(-\int_r^t D(s)ds - \int_0^{t-r} \beta(s)ds\right) \varphi(t, t-r) dr = \int_0^t x(r) \varphi(r, 0) dr$$

It follows from the above equations and (A1) that equation (4.2) holds.

Suppose that $\bar{f} \in C^0\left([0,T] \times \mathbb{R}_+\right) \cap L^1\left([0,T] \times \mathbb{R}_+\right)$ with $\sup_{t \in [0,T]}\left(\|\bar{f}[t]\|_1\right) < +\infty$ satisfies equation (4.2) for all $\varphi \in C^1\left([0,T] \times \mathbb{R}_+\right) \cap L^\infty\left([0,T] \times \mathbb{R}_+\right)$ with $\left(\frac{\partial \varphi}{\partial a} + \frac{\partial \varphi}{\partial t}\right) \in L^\infty\left([0,T] \times \mathbb{R}_+\right)$ and $t \in [0,T]$. We next show that $f \equiv \bar{f}$. Define the function $u = f - \bar{f}$. Then $u \in C^0\left([0,T] \times \mathbb{R}_+\right)$ with $\sup_{t \in [0,T]}\left(\|u[t]\|_1\right) < +\infty$ and $u$ satisfies the following equation for all $\varphi \in C^1\left([0,T] \times \mathbb{R}_+\right) \cap L^\infty\left([0,T] \times \mathbb{R}_+\right)$ with $\left(\frac{\partial \varphi}{\partial a} + \frac{\partial \varphi}{\partial t}\right) \in L^\infty\left([0,T] \times \mathbb{R}_+\right)$ and $t \in [0,T]$:

$$\int_0^t \int_0^{+\infty} \left((\beta(a) + D(s))\varphi(s,a) - \frac{\partial \varphi}{\partial a}(s,a) - \frac{\partial \varphi}{\partial s}(s,a)\right) u(s,a) da\,ds \quad (A5)$$
$$+ \int_0^{+\infty} u(t,a) \varphi(t,a) da = 0$$

Let arbitrary $\tau \in [0,T]$, $\tilde{D} \in C^0\left([0,T]\right)$ and $g \in C^1\left([-T, +\infty)\right) \cap L^\infty\left([-T, +\infty)\right)$ be given. Define the function $\varphi \in C^1\left([0,T] \times \mathbb{R}_+\right) \cap L^\infty\left([0,T] \times \mathbb{R}_+\right)$:

$$\varphi(t,a) = \exp\left(-\int_t^\tau \tilde{D}(s)ds - \int_a^{a+\tau-t} \beta(s)ds\right) g(a+\tau-t) \quad (A6)$$

We notice that the function $\varphi$ defined by (A6) satisfies $\frac{\partial \varphi}{\partial a}(t,a) + \frac{\partial \varphi}{\partial t}(t,a) = \left(\tilde{D}(t) + \beta(a)\right)\varphi(t,a)$ for all $t \in [0,T]$, $a \geq 0$. Since $\tilde{D} \in C^0([0,T])$, $\beta \in C^0(\mathbb{R}_+) \cap L^\infty(\mathbb{R}_+)$ and $\varphi \in L^\infty([0,T] \times \mathbb{R}_+)$, we obtain that $\left(\frac{\partial \varphi}{\partial a} + \frac{\partial \varphi}{\partial t}\right) \in L^\infty([0,T] \times \mathbb{R}_+)$. Thus, we get from (A5) and (A6) for $t \in [0,T]$:

$$\int_0^t \left(D(s) - \tilde{D}(s)\right) \int_0^{+\infty} u(s,a) \exp\left(-\int_s^\tau \tilde{D}(l)dl - \int_a^{a+\tau-s} \beta(l)dl\right) g(a+\tau-s) da\, ds$$
$$+ \int_0^{+\infty} u(t,a) \exp\left(-\int_t^\tau \tilde{D}(s)ds - \int_a^{a+\tau-t} \beta(s)ds\right) g(a+\tau-t) da = 0 \quad (A7)$$

Setting $t = \tau$, we get from (A7):

$$\int_0^\tau \left(D(s) - \tilde{D}(s)\right) \int_0^{+\infty} u(s,a) \exp\left(-\int_s^\tau \tilde{D}(l)dl - \int_a^{a+\tau-s} \beta(l)dl\right) g(a+\tau-s) da\, ds$$
$$+ \int_0^{+\infty} u(t,a) g(a) da = 0 \quad (A8)$$

Using the fact that for every $\varepsilon > 0$ there exists $\tilde{D} \in C^0([0,T])$ with $\int_0^T |D(s) - \tilde{D}(s)| ds < \varepsilon$ and $\|\tilde{D}\|_\infty \leq \|D\|_\infty$, the fact that $\sup_{t \in [0,T]} \left(\|u[t]\|_1\right) < +\infty$ and the fact that $\tilde{D} \in C^0([0,T])$ is arbitrary, we obtain from (A8):

$$\int_0^{+\infty} u(\tau,a) g(a) da = 0 \quad (A9)$$

Using the fact that $g \in C^1([-T,+\infty)) \cap L^\infty([-T,+\infty))$ is arbitrary and exploiting the fact that $u[\tau] \in C^0(\mathbb{R}_+)$ with $\|u[\tau]\|_1 < +\infty$, we establish from (A9) and Corollary 4.24 on page 110 in [5] that $u[\tau] = 0$. Since $\tau \in [0,T]$ is arbitrary, we conclude that $u(t,a) = 0$ for all $t \in [0,T]$, $a \geq 0$.

The proof is complete. ◁

**Proof of Lemma 2:** Since $(f,S):[0,T] \to X$ is a continuous mapping with metric given by (2.5), it follows that the mapping $[0,T] \ni t \to \int_0^{+\infty} (\beta(a)+D(t))f(t,a)da$ is of class $L^\infty([0,T])$. Therefore, equation (2.7) with $\varphi(t,a) \equiv 1$ implies that the mapping $[0,T] \ni t \to \int_0^{+\infty} f(t,a)da$ is absolutely continuous and satisfies for $t \in [0,T]$ a.e.:

$$\frac{d}{dt}\left(\int_0^{+\infty} f(t,a)da\right) = f(t,0) - \int_0^{+\infty} (\beta(a)+D(t))f(t,a)da \qquad (A10)$$

Since $f(t,a) > 0$ for all $a \geq 0$, it holds that

$$\|f[t]\|_1 = \int_0^{+\infty} f(t,a)da \qquad (A11)$$

Combining (A10), (A11) and (2.2), we get for $t \in [0,T]$ a.e.:

$$\frac{d}{dt}\left(\|f[t]\|_1\right) = \mu(S(t))\int_0^{+\infty} k(a)f(t,a)da - D(t)\|f[t]\|_1 - \int_0^{+\infty} \beta(a)f(t,a)da \qquad (A12)$$

Exploiting the facts that $D(t) \geq 0$ for $t \geq 0$ a.e., $\beta(a) \geq 0$ for all $a \geq 0$, we get for $t \in [0,T]$ a.e.:

$$\frac{d}{dt}\left(\|f[t]\|_1\right) \leq M\|k\|_\infty \|f[t]\|_1 \qquad (A13)$$

where $M = \max_{S \in [0,S_{in}]}(\mu(S))$. The differential inequality (A13) implies estimate (4.3).

Using (2.3) and the facts that $\mu(S) \geq 0$ for all $S \geq 0$, $q(a) \geq 0$ for all $a \geq 0$, we obtain the following differential inequality for $t \in [0,T]$ a.e.:

$$\dot{S}(t) \leq D(t)(S_{in} - S(t)) \qquad (A14)$$

The differential inequality (A14) implies the second estimate (4.4).

Using (4.3), (2.3) and the facts that $\mu(S) \leq \Gamma S$ holds for all $S \in [0, S_{in}]$, $\mu(S) \geq 0$ for all $S \geq 0$, $q(a) \geq 0$ for all $a \geq 0$, we obtain the following differential inequalities for $t \in [0,T]$ a.e.:

$$\dot{S}(t) \geq D(t)\left(S_{in} - S(t)\right) - \mu(S(t))\|q\|_\infty \int_0^{+\infty} f(t,a)da$$

$$\geq D(t)\left(S_{in} - S(t)\right) - \mu(S(t))\|q\|_\infty \exp\left(M\|k\|_\infty t\right)\|f_0\|_1 \qquad (A15)$$

$$\geq D(t)\left(S_{in} - S(t)\right) - \Gamma\|q\|_\infty \exp\left(M\|k\|_\infty t\right)\|f_0\|_1 S(t)$$

The differential inequality (A15) implies the following estimates for all $t \in [0,T]$:

$$S(t) \geq \exp\left(-\int_0^t D(s)ds - \Gamma\|q\|_\infty \|f_0\|_1 \frac{\exp(M\|k\|_\infty t) - 1}{M\|k\|_\infty}\right) S_0$$

$$+ S_{in} \int_0^t \exp\left(-\int_\tau^t D(s)ds - \Gamma\|q\|_\infty \|f_0\|_1 \frac{\exp(M\|k\|_\infty t) - \exp(M\|k\|_\infty \tau)}{M\|k\|_\infty}\right) D(\tau)d\tau$$

$$\geq \exp\left(-\Gamma\|q\|_\infty \|f_0\|_1 \frac{\exp(M\|k\|_\infty t) - 1}{M\|k\|_\infty}\right) \left(\exp\left(-\int_0^t D(s)ds\right) S_0 + S_{in}\left(1 - \exp\left(-\int_0^t D(s)ds\right)\right)\right)$$

Since $\exp\left(-\int_0^t D(s)ds\right) S_0 + S_{in}\left(1 - \exp\left(-\int_0^t D(s)ds\right)\right) \geq S_0$ (a consequence of the facts that $D(t) \geq 0$ for $t \geq 0$ a.e. and $S_0 \in (0, S_{in})$), the above estimate implies the first estimate (4.4). The proof is complete. ◁